\documentclass[oneside,a4paper,12pt,reqno]{amsart}
\usepackage{upref,amsfonts,amsxtra,amssymb}

\newcommand{\Z}{\mathbb{Z}}

\newcommand{\R}{\mathbb{R}}
\newtheorem{theorem}{Theorem}[section]
\newtheorem{lemma}[theorem]{Lemma}

\newtheorem{cor}[theorem]{Corollary}
\title{Balanced words in higher dimensions}
\author{Siddhartha Bhattacharya}
\address{School of Mathhematics, Tata Institute of Fundamental
  Research, Mumbai 400005, India}
\email{siddhart@math.tifr.res.in}
\subjclass[2010]{5B99}
\date{}
\begin{document}
\maketitle
\begin{abstract}
For $d\ge 1$, a word $w\in \{ 0,1\}^{\Z^d}$ is called balanced if
there exists $M > 0$ such that for any two rectangles $R,
R^{'}\subset\Z^d$ that are translates of each other, the number of
occurrences of the symbol $1$ in $R$ and $R^{'}$ differ by at most
$M$.  It is known that for every balanced word $w$, the asymptotic frequency of
the symbol $1$ ( called the density of $w$ ) exists. In this paper we
show that there exist two dimensional balanced words with
irrational densities, answering a question raised by Berth\'e and Tijdeman. 
\end{abstract}

\section{Introduction}
A word $w\in \{ 0,1\}^{\mathbb Z}$ is called {\it $k$-balanced} if for
any two blocks $B, B^{'}\subset\Z$ of equal length, the number of
occurrences of the symbol $1$ in $B$ and $B^{'}$ differ by at most $k$.
We will call a word $w$ {\it balanced} if it is $k$-balanced for some
$k\ge 1$. In the literature $1$-balanced words are often called
balanced. However, following \cite{BT}, we will reserve this term for the
weaker property. 

Balanced words occur naturally in many different areas, including ergodic
theory, number theory, and theoretical computer science (\cite{GJWZ},
\cite{GZ},\cite{T}). 
The study of balanced words was initiated by Morse and Hedlund
(\cite{MH1}, \cite{MH2}). They obtained a classification theorem for
  $1$-balanced words involving the density (asymptotic
  frequency of the symbol 1) of the underlying
  word. In particular, they proved that any $1$-balanced word 
$w\in \{ 0,1\}^{\Z }$ admits a density, and $1$-balanced
  words in $\{ 0,1\}^{\Z}$ with irrational densities 
are Sturmian words corresponding to certain codings of irrational
rotations of the circle.

In \cite{BT}, Berth\'e and Tijdeman studied balanced words in 
$\{ 0,1\}^{\Z^d}$ for $d >1$. They showed that 
higher dimensional balanced words also admit densities, but 
for $d\ge 2$, the density of a $1$-balanced word is always rational
(\cite[the main corollary]{BT}). The following question has been raised by
several authors (\cite{BT},\cite[Conjecture 2]{V}) :
  
 Question : {\it For $d\ge 2$, does there exist a balanced word 
in $\{ 0,1\}^{\Z^d}$ with irrational  density ? }

In this paper we give an affirmative answer to this
question. 
\begin{theorem}
For any $a\in[0,1]$ there exists a balanced word $w$ in $\{0,1\}^{\Z^2}$
with density $a$.
\end{theorem}

If $a$ is an irrational number sufficiently close
to ${\frac{1}{2}}$ then our proof shows that 
there exists a $k$-balanced word $w$ in 
$\{ 0,1\}^{\Z^2}$ with $k = 33$ that has density $a$.
 
\section{Density of balanced words}
For a set $B$, $|B|$ will denote the cardinality of $B$. For
$a,b\in\Z$ with $a\le b$ we will denote the set $\{ a, a+1, \ldots ,
b-1\}$ by $[a, b)$. A {\it
  rectangle\/} $R$ in 
$\Z^d$ is a set of the form $[a_1, b_1)\times
\cdots \times [a_d, b_d)$. 
For $d\ge 1$, $W_d$ will denote the set $\{ 0,1\}^{\Z^d}$, the
collection of all $d$-dimensional words with $\{ 0,1\}$ as the
alphabet set.
If $w\in W_d$, we define $S_w \subset \Z^d$ by 
$$S_w = \{ g\in\Z^d : w(g) = 1 \}.$$
 A word $w$ is said to have a
{\it density \/} if there exists $a\in [0,1]$ such that 
$|S_{w}\cap R_n|/|R_n|\mapsto a$
as $n\mapsto\infty$
for every
sequence of rectangles $\{ R_n\}$ with $|R_n|\mapsto \infty$.

Suppose $d\ge 1$, $a\in [0,1]$ and $f$ is a function from 
$\Z^d$ to the two element set $\{a, a-1\}$. We define $w(f)\in W_d$
by $w(f)(g) = a - f(g)$.
\begin{lemma}\label{function}
Let $a$,$f$ and $w(f)$ be as above. Suppose there exists $M >0$ such that
for every rectangle $R$,
$$|\sum_{g\in R}f(g)| \le M.$$
Then $w(f)$ is a balanced word with density $a$.   
\end{lemma}  
{\bf Proof.} Since $f(g) = a - w(f)(g)$ for all $g\in \Z^d$, it follows
that for any rectangle $R$,
$$\sum_{g\in  R}f(g) = a |R| - |S_{w(f)}\cap R|.$$
Hence $|{\frac{|S_{w(f)}\cap R|}{ |R|}} -  a| \le {\frac{M}{|R|}}$.
 This implies that
$|S_{w(f)}\cap R|/|R|\mapsto
a$ as $|R|\mapsto \infty$, i.e., $w(f)$ has density $ a$.
If $R^{'}$ is a translate of $R$ then the above equation 
also shows that 
$$  ||S_{w(f)}\cap R| -  |S_{w(f)}\cap R^{'}||  = |\sum_{g\in  R}f(g) - 
\sum_{g\in  R^{'}}f(g)| \le 2M.$$
Hence $w(f)$ is a $2M$-balanced word. $\Box$

\medskip
For $d\ge 1$, let $Q_d$ denote the ring ${\mathbb
  Z}[u_{1},\ldots , u_{d}]$, the polynomial ring 
in $d$ commuting variables with integer
coefficients. For $g = (n_1 ,\ldots ,n_d)\in {\mathbb Z}^d$, we
will denote the monomial
$u_{1}^{n_{1}}\cdots u_{d}^{n_{d}}$ by $u^{g}$. Every
element of $Q_d$ can be expressed as $\sum c_{g} u^{g}$,
with $c_{g}\in \Z$ and $c_{g} = 0$ for all but
finitely many $g$. For $d\ge 1$, $l^{\infty}(\Z^d)$ will denote the
space of all bounded functions from $\Z^d$ to $\R$.   
For $p = \sum_{g} c_{g}u^{g}$ and $f\in
l^{\infty}(\Z^d)$ we define $p\cdot f\in l^{\infty}(\Z^d)$ by 
$$ p\cdot f(x) = \sum_{ g} c_{g} f(x + g).$$  

\begin{theorem}\label{unip}
Suppose $d\ge 1$ and $a\in [0,1]$. If there exists $f\in l^{\infty}(\Z)$ 
 such that 
$((u - 1)^d\cdot f )(m) \in \{ a, a-1\}$ for all $m\in \Z$, then
there exists a balanced word $w\in W_d$  with density $a$. 
\end{theorem}
{\it Proof. \/} Let $\theta : \Z^d\rightarrow \Z$ denote the homomorphism
 defined by 
$$\theta((n_1, \ldots , n_d) = \sum_{1}^{d}
n_i,$$ 
and let
$\theta^{*}$ denote the ring homomorphism from $Q_d$ to $Q_1$ induced
by $\theta$, i.e., $\theta^{*}(\sum c_gu^{g}) = \sum c_gu^{\theta(g)}$.
We note that $f\circ \theta$ is an element of $l^{\infty}(\Z^d)$, and 
for any $g,g^{'}\in\Z^d$,
$$u^{g}\cdot (f\circ \theta)(g^{'}) = f\circ\theta(g + g^{'}) =
u^{\theta(g)}f(\theta(g^{'})).$$ 
Hence $p\cdot (f\circ \theta) = (\theta^{*}(p) \cdot f)\circ \theta$
for all $p\in Q_d$. We define $h\in l^{\infty}(\Z^d)$ and $p_0\in Q_d$
by
$$ p_0 = \prod_{i = 1}^{d} (u_i - 1),\ \ h = p_0\cdot (f\circ
\theta).$$
It is easy to see that $\theta^{*}(p_0) = (u -1)^d$. Since the image
of $\theta^{*}(p_0) \cdot f$ is contained in $\{ a, a-1\}$, we deduce
that the image of $h = (\theta^{*}(p_0) \cdot f)\circ \theta$ is also 
contained in $\{ a, a-1\}$. Let $R = \prod_{i} [a_i, b_i)$ be an
arbitrary rectangle in $\Z^d$. Let $q\in Q_d$ denote the polynomial 
$\sum_{g\in R}u^{g}$. Then 
$$\sum_{g\in R} h(g) = \sum_{g\in R} (u^{g}\cdot h)(0) = (q\cdot
h)(0).$$
 We note that $q = q_1\cdots q_d$, where
$ q_i =  \sum_{k\in [a_i, b_i)} u_{i}^{k}$.  As $q\cdot h = qp_0\cdot
(f\circ \theta)$, this implies that 
$$ q\cdot h = (\prod_{i=1}^{d} (u_{i}^{b_i} - u_{i}^{a_i}))\cdot(f\circ \theta).$$
We note that $||s||_{\infty} = ||u^{g}s||_{\infty}$ 
for any $s\in l^{\infty}(\Z^d)$ and $g\in\Z^d$. Since the right hand
side of the above equality has $2^d$ terms of the form
$\pm u^{g}\cdot (f\circ \theta)$, we conclude that 
$||q\cdot h||_{\infty} \le 2^d||f||_{\infty}$. In particular, 
$\sum_{g\in R} h(g) = (q\cdot h)(0)\le 2^d||f||_{\infty}$. As
$R\subset\Z^d$ is arbitrary, the given assertion 
follows from the previous lemma. $\Box$

\medskip
\noindent
It is easy to see that for any $f\in l^{\infty}(\Z)$ 
the value of  $(u - 1)^d\cdot f$ at $m$
depends only on $\{ f(m), f(m+1),\ldots \}$. Hence if
${\mathcal S}$ denotes the space of all sequences taking values in
$\R$ then $x\mapsto (u - 1)^d\cdot x $ is a well defined map from 
${\mathcal S}$ to ${\mathcal S}$. We note the following consequence of
the previous result :
\begin{cor}\label{direction}
Suppose $d\ge 1$ and $a\in [0,1]$. If there exists a bounded sequence $x$ 
such that $x_1 = \cdots = x_{d} = 0$ and $((u - 1)^d\cdot
x)_n$ lies in $\{ a, a-1\}$ for all $n\ge 1$, then 
there exists a balanced word $w\in W_d$  with density $a$.
\end{cor} 
{\it Proof.\/} In view of the previous theorem it is enough to
construct $f\in l^{\infty}(\Z)$ such that 
$((u - 1)^d\cdot f )(m) \in \{ a, a-1\}$ for all $m\in \Z$.
We note that $(u-1)^d = \sum c_iu^i$, where 
$c_i = (-1)^{d-i}{\binom{d}{i}}$. Suppose $d$ is even. 
 We define $f : \Z\rightarrow \R$
by $f(m) = x_m$ if $m\ge 1$, and $f(m) = x_{d + 1-m}$ if $m\le 0$. Clearly
$f$ is an element of $l^{\infty}(\Z)$. 
Since $((u - 1)^d\cdot
x)_n$ lies in $\{ a, a-1\}$ for all $n\ge 1$, and $c_i =
c_{d-i}$ for all $i$, it follows that 
$((u - 1)^d\cdot f )(m)$
lies in $\{ a, a-1\}$ for all $m\in \Z$. 
If $d$ is odd, we define $f(m) = x_m$ if $m\ge 1$, and 
$f(m) = - x_{d + 1-m}$ if $m\le 0$. Since $c_i =
- c_{d-i}$ for all $i$, we deduce that 
$((u - 1)^d\cdot f )(m)$
lies in $\{ a, a-1\}$ for all $m\in \Z$. $\Box$ 

\section{Two dimensional balanced words}
Now we consider the case when $d = 2$. 
Let $h$ be an arbitrary
function from $\R^2$ to $\{ a, a -1\}$. We define a function 
$T_h : \R^2\rightarrow\R^2$ by
$$T_{h}(x,y) = (x+ y, y + h(x,y)).$$
We define $(x_n, y_n) = T_{h}^{n}(0,0)$, $x = \{ x_i\}$ and $y = \{
y_i\}$. Since $x_{n+1} =
x_n + y_n$, it follows that $(u-1)\cdot x= y$. Hence for all $n\ge 1$,
$$ ((u-1)^{2}\cdot x)_n = y_{n+1} - y_n = h(x_n,
y_n)\in \{ a, a-1\}.$$
As $x_1 = x_2 = 0$,   
 Theorem 1.1 follows from Corollary
\ref{direction} and the following result :
 
\begin{theorem}
For any $a\in [0,1]$, there exists $h : \R^2\rightarrow \{ a, a-1\}$
such that the set $\{ T_{h}^{n}(0,0) : n\ge 1\}\subset \R^2$ is bounded. 
\end{theorem}
{\it Proof. \/} 
We define $h : \R^2\rightarrow \{ a, a-1\}$ as follows :   
$h(x,y) = a -1$ if $y+a > 1$, or
if both $x$ and $y + a$ are positive. Otherwise, $h(x,y) = a$.
We define $(x_n, y_n) = T_{h}^{n}(0,0)$ and claim that $-1\le y_i\le 1$
for all $i$.
Suppose this is not the case. If there
exists $m$ such that $y_m > 1$, we set $l = \min \{ i : y_i > 1\}$.
Then $y_{l} = y_{l-1} + h(x_{l-1}, y_{l-1}) > 1$. From the minimality
of $l$ we deduce that $h(x_{l-1}, y_{l-1}) > 0$. As the range of $h$
is $\{ a, a-1\}$, this shows that $h(x_{l-1}, y_{l-1}) = a$. But this
can happen only if $y_{l-1} + a \le 1$, i.e., $y_l \le 1$. Since this
contradicts $y_{l} > 1$, we conclude that  $y_i \le 1$ for all $i$.
Now suppose $y_i < -1$ for some $i$. 
We define $l = \min \{ i : y_i < - 1\}$. Clearly $y_{l-1} > y_l$,
i.e., $h(x_{l-1}, y_{l-1}) = y_l - y_{l-1}$ is negative. Hence 
$h(x_{l-1}, y_{l-1}) = a -1$
and $y_{l-1} = y_l + 1 - a$. But this implies that 
$$y_{l-1} + a  = y_l + 1 < 0,$$
which is not possible since by construction $h(x_{l-1}, y_{l-1}) = a$ whenever 
$ y_{l-1} + a \le 0$. This proves the claim.

Now we will show that $x_i \le \alpha = {\frac{1}{1-a}} + 2 $
for all $i$. Suppose
this is not the case. We pick any $k$ such that $x_k > \alpha$ and
define 
$$l =  \max \{ i \le k : x_i \le 0 \}.$$
Clearly $x_{i} > 0$ for
$i\in [l+1, k]$. We look at the finite sequence
$y_l, y_{l+1}, \ldots, y_{k-1}$.
If $y_i > 0$ for some $i\in [l,k-1]$, then $y_i + a$ is also positive.
Since $x_i > 0$, this implies that $h(x_i, y_i) = a -1$. Hence 
$y_{i+1} = y_i + a - 1$. If $y_i \le 0$ then there are two possibilities.
In the first case, when $y_i + a > 0$, from the definition of $h$ it
follows that $ h(x_i, y_i) = a - 1$. In particular, $y_{i+1} \le
y_i\le 0$. In the second case, when  $y_i + a \le 0$, it is easy to
see that $y_{i+1} = y_i + h(x_i, y_i)\le 0$ irrespective of whether 
$h(x_i, y_i) = a$ or $a-1$. Combining these two cases, we observe that 
$y_{i+1} \le 0$ whenever $y_i\le 0$. Therefore the sequence 
$y_l, y_{l+1}, \ldots, y_{k-1}$ decreases by $1 -a$ at each step as
long as $y_i$ is positive, and once it becomes non-positive it remains
 non-positive.
So there exists a unique $j\in
[l,k-1]$
such that $ y_i  > 0$ for $ l\le i < j$ and $y_{i}\le 0$ for
$i\ge j$. 
Since $x_{i+1} - x_i = y_i$, it follows that
$$x_{j+ 1} \ge x_{j+ 2} \ge \cdots \ge x_{k} > \alpha.$$
We note that $y_l\le 1$, $y_{j-1} > 0$, and $y_{i+1} = y_i - (1-a)$ for
$i< j$. This shows that $(j-1-l)\le {\frac{1}{1 -a}}$. 
From the previous claim  we see that $x_{i+1} = x_i + y_i \le x_i+1$
for all $i$. As $x_l \le 0$, we obtain that 
$$\alpha <  x_{j+1} \le j+1-l\le 2 +  {\frac{1}{1 -a}}$$.
This contradiction shows that $x_i\le \alpha$ for all $i$. 

To complete the proof of the given assertion we need to show that the
sequence $\{ x_i\}$ is also bounded from below. We set 
$\beta = 2 + {\frac{1}{a}}$ and claim that $x_i \ge - \beta$ for all
$i$. Suppose this is not the case. 
We pick $k$ such that $x_k < - \beta$ and
define 
$$l =  \max \{ i \le k : x_i \ge 0 \}.$$
Clearly $x_{i} < 0$ for
$i\in [l+1, k]$. As before, we look at the finite sequence
$y_l, y_{l+1}, \ldots, y_{k-1}$.
If $y_i < 0$ for some $i\in [l,k-1]$, then $y_i + a < 1$.
Since $x_i \le 0$, this implies that $h(x_i, y_i) = a$. Hence 
$y_{i+1} = y_i + a$. If $y_i \ge 0$ then there are two possibilities.
In the first case, when $y_i + a \le 1$, from the definition of $h$ it
follows that $ h(x_i, y_i) = a$. In particular, $y_{i+1} >
y_i\ge 0$. 
In the second case, when  $y_i + a > 1$, it is easy to
see that $y_{i+1} = y_i + h(x_i, y_i)> 0$ irrespective of whether 
$h(x_i, y_i) = a$ or $a-1$. Combining these two cases, we observe that 
$y_{i+1} \ge 0$ whenever $y_i\ge 0$. Therefore the sequence 
$y_l, y_{l+1}, \ldots, y_{k-1}$ increases by $a$ at each step as
long as $y_i$ is negative, and once it becomes non-negative it remains
 non-negative.
So there exists a unique $j\in
[l,k-1]$
such that $ y_i  < 0$ for $ l\le i < j$ and $y_{i}\ge 0$ for
$i\ge j$. 
Since $x_{i+1} - x_i = y_i$, it follows that
$$x_{j+ 1} \le x_{j+ 2} \le \cdots \le x_{k} < -\beta.$$
We note that $y_l\ge -1$, $y_{j-1} < 0$, and $y_{i+1} = y_i + a$ for
$i< j$. This shows that $(j-1-l)\le {\frac{1}{a}}$. 
From the previous claim  we see that $x_{i+1} = x_i + y_i \ge x_i- 1$
for all $i$. As $x_l \ge 0$, we obtain that 
$$-\beta > x_{j+1} \ge l- 1 - j\ge - (2+{\frac{1}{a}}).$$
This contradiction shows that $x_i\ge -\beta$ for all $i$. $\Box$

\end{document}